\documentclass{article}
\usepackage{amsthm}
\usepackage{amsfonts}

\begin{document}

\centerline{\bf Some properties of $G$-schemes.}
\bigskip
\bigskip
\bigskip

Tsemo Aristide

College Boreal

1 Yonge Street, Toronto

tsemo58@yahoo.ca

\bigskip
\bigskip
\bigskip

\centerline{\bf Abstract.}

\bigskip
\bigskip

In this paper, we continue to adapt the theories of spectra and schemes developed by Grothendieck in algebraic geometry to the category of groups. Let $G$ be a group,  and $(H,f_G^H)$ and object of the comma category $C(G)$. In [5], we have defined on the set $Spec_G(H)$ of prime ideals of $(H,f_G^H)$ a topology. In this paper, we define another notion of prime ideals to which we associate a spectrum endowed with a topology. We study some properties and objects associated to  these spectra; amongst them we can quote, irreducibility, the radical, the structural sheaf, $G$-varieties and $G$-schemes.
\bigskip
\bigskip

{\bf 1. Introduction.}

\bigskip
\bigskip

 Let $G$ be a group and $F[X_i, i\in I]$ the free group generated by the set $(X_i)_{i\in I}$. We denote by $G[X_i, i\in I]$ the free product of  $G$ and $F[X_i, i\in I]$. A $G$-group is defined by a group $H$, and a morphism of groups $f_H:G\rightarrow H$. An element $f$ of $G[X_1,...,X_n]$ defines a  function  $e_f:H^n\rightarrow H$ such that $f(h_1,..,h_n)$ is obtained by substituting $X_i$ by $h_i$ in the expression of $f$. 
  We say that $f$ is a polynomial function defined on $H^n$. Let $S$ be a subset of $G[X_1,..,X_n]$,   the set of zeroes of $S$,  or equivalently the algebraic subset defined by $S$ is the subset $V_H(S)$ of $H^n$ defined by $V_H(S)=\{(h_1,..,h_n)\in H^n: f(h_1,..,h_n)=1, f\in S\}$. In [1], the authors 
     have defined a (Zariski) topology on $H^n$ whose set of closed subsets is generated  by the subbasis of algebraic subsets. They have developed in this context notions similar to the notions defined  in classical algebraic geometry theory and commutative algebra like the notion of prime ideals, that we will call here $T_1$-prime ideals here.
    
     In [5], we have started to adapt the ideas of Grothendieck's scheme in this context. Consider a $G$-group $f_H:G\rightarrow H$ such that $f_H$ is injective.  We have endowed  $Spec^1_G(H)$, the set of $T_1$-prime ideals  of $H$ with a  topology for which, a subset $C\subset Spec^1_G(H)$ is closed if and only if, there exists a normal subgroup $I$ of $H$ such that $C=V^1_H(I)=\{P\in Spec^1_G(H), I\subset P\}$. The space $Spec^1_G(H)$ endowed with this topology is called a  $T_1$-$G$-affine spectrum.
     
      In the second part of this paper, we do not assume that $f_G^H$ is injective, and we define another notion of prime ideals that we call $T_2$-prime ideals and denote by $Spec_G^2(H)$ the set of $T_2$-prime ideals. 
      We define the notions of  $T_i$-topologies $i=1,2$. We also show that the set of $T_1$-divisors of zero of $G[X_i,i\in I]$ is contained in $G$  if $G$ is different of the trivial group and $\mathbb{Z}/2$. We also study and characterize irreducible $T_i$-$G$-spectra. In the third part, we study affine $G$-varieties that we endow with a Zariski topology when $G$ is an integral domain. In the forth part we study the structural sheaf of an affine $T_i$-$G$-scheme $Spec^i_G(H), i=1,2$, we describe its stalks and its global sections if the intersection of the irreducible components of $Spec^i_G(H)$ is not empty. In the last part, we study  $G$-schemes.

\bigskip
\bigskip

{\bf Plan.}

\bigskip

1. Introduction.

2. Some properties of groups and affine $G$-spectra.

3. $G$-varieties.

4. Structural preshaf and structural sheaf.

5. $G$-schemes.    
     
\bigskip
\bigskip

{\bf 2. Some properties of groups and affine $G$-spectra.}
\bigskip

{\bf 2.1 Definitions and  Notations.}

\bigskip

Let $Group$ be the category of groups and $H$ an object of $Group$. Let $L$ and $L'$ be two subgroups of $H$, we denote by $[L,L']$ the subgroup of $H$ generated by $\{ll'l^{-1}{l'}^{-1}, l\in L, l'\in L'\}$.

For any group $G$, we denote by $C(G)$ the comma category of $G$. Recall that an object $(H,f_G^H)$   of $C(G)$ is a morphism of groups $f_G^H: G\rightarrow H$, and a morphism $m:(H,f_G^H)\rightarrow (H',f^{H'}_G)$ is a morphism of groups $m:H\rightarrow H'$ such that $m\circ f_G^H=f^{H'}_G$. 

\bigskip

Let $(H,f_G^H)$ be an object of $C(G)$, and $x$ an element of $H$, we denote by $G(x)$ the subgroup of $H$ generated by $\{f_G^H(g)x{f_G^H(g)}^{-1}, g\in G\}$.

An element $x\in H$  is a $T_1$-divisor of zero if and only if,  $x$ is distinct of $1$, and there exists an element $y\in H$  distinct of $1$, such that  $[G(x),G(y)]=1$.  We say that $x$ is a $T_2$-divisor of zero if and only if there exists an element $y$ distinct of $1$ such that $G(x)\cap G(y)=1$.

An element $x$ of $(H,f_G^H)$ is nilpotent if and only if  $G(x)$ is nilpotent. We denote by $Nil_G(H)$ the set of nilpotent elements of $H$.

Remark that if $x$ is a nilpotent element distinct of $1$, and $y$ is an element of the center of $G(x)$ distinct of $1$,  $[G(x),G(y)]=1$. We deduce that $x$ is a $T_1$-divisor of zero.

An ideal of $(H,f_G^H)$ is a normal subgroup $I$ of $H$,  distinct of $H$.
Remark that if $I$ is an ideal of $(H,f_G^H)$, the composition of $f_G^H$ and the quotient morphism $p_I:H\rightarrow H/I$ defines an object $(H/I,p_I\circ f_G^H)$ of $C(G)$.

An object $(H,f_H)$ is $T_i$-integral if it does not have a $T_i$-divisor of zero, $i=1,2$.

An ideal $I$ of $(H,f_G^H)$ is a $T_i$-prime ideal if and only if $(H/I,p_I\circ f_G^H)$ is $T_i$-integral   $i=1,2$. Remark that the ideal $I$ is $T_1$-prime if and only if for every $x,y\in H$, $[G(x),G(y)]\subset I$ implies that $x\in I$ or $y\in I$. The ideal $I$ is $T_2$-prime if and only if for every $x,y\in H$, $G(x)\cap G(y)\subset I$ implies that $x\in I$ or $y\in I$. Remark that the fact that $I$ is a $T_2$-prime ideal implies that $I$ is a $T_1$-prime ideal.

The  $T_i$-$G$-spectrum  $Spec^i_G(H)$ is the set of $T_i$-prime ideals of $(H,f_G^H)$ , $i=1,2$. Let $P$ be an element of $Spec^i_G(H)$, we will call the quotient $H/P$, the local group of $P$. Let $I$ be a normal subgroup of $(H,f_G^H)$, we denote  by $V^i_H(I)$ the set of $T_i$-prime ideals which contain $I$.  Remark that
 the spectrum $Spec^i_G\{1\}$ of the  object  $(\{1\},f_G^1)$  is empty.

\bigskip

{\bf Theorem 2.1.}

{\it Suppose that  $G$ is distinct of the trivial group and the cyclic group $\mathbb{Z}/2$. If an element $x\in G[X_i,i\in I]$ is  a $T_1$-divisor of zero then $x\in G$.}

\begin{proof}
Every element of $G[X_i,i\in I]$ distinct of the identity is a divisor of zero if $G=\{1\}$. Suppose that $G=\mathbb{Z}/2$, and denote by $g$ the element of $\mathbb{Z}/2$ distinct of the identity. Let $i\in I$, $G(gX_igX_i^{-1})$ is the group generated by $u=gX_igX_i^{-1}$ and its conjugate by $g$, $v=X_igX_i^{-1}g$. Since $u=v^{-1}$, $G(gX_igX_i^{-1})$ is a cyclic group. We deduce that $[G(gX_igX_i^{-1}),G(gX_igX_i^{-1})]=1$. This implies that $gX_igX_i^{-1}$ is a divisor of zero of $G[X_i, i\in I]$.

Suppose that $G$ is neither the trivial group, nor $\mathbb{Z}/2$. Let $x$ be a divisor of zero of $G[X_i,i\in I]$.
Consider an element $y$ distinct of $1$ such that $[G(x),G(y)]=1$. The Corollary 4.1.6 of [4] implies that either $x$ and $y$ are contained in the same conjugate of $G$, or there exists $u$ such that $x=u^n, y=u^m$.

Suppose that $x$ and $y$ are in the same conjugate of $G$. There exists $w$ such that $x=wgw^{-1}, y=wg'w^{-1}, g,g'\in G$. 
  We can write the reduced form of $w$:
  
  1. $w=X_{i_1}^{l_1}g_{i_2}X_{i_2}^{l_2}..X_{i_n}^{l_n}$, or

2. $w= g_{i_1}X_{i_1}^{l_1}g_{i_2}X_{i_2}^{l_2}..X_{i_n}^{l_n},$ or

3.  $w=X_{i_1}^{l_1}g_{i_2}X_{i_2}^{l_2}..X_{i_n}^{l_n}g_{i_{n+1}}$, or

4.   $w=g_{i_1}X_{i_1}^{l_1}g_{i_2}X_{i_2}^{l_2}..X_{i_n}^{l_n}g_{i_{n+1}}$.
   
  For the first  and the third cases, the reduced forms of $x$, and $y$, have at their ends elements of the form $X_i$ or $X_i^{-1}$. For every $h\neq 1\in G$, $[hxh^{-1},y]\neq 1$. This contradicts the fact that $[G(x),G(y)]\neq 1$.
  
  For the second and the forth cases, the reduced forms of $g_{i_1}^{-1}xg_{i_1}$ and $g_{i_1}^{-1}yg_{i_1}$ have at their ends elements of the form $X_i$ and $X_i^{-1}$, for every $h\neq 1\in G$, $[hg_{i_1}^{-1}xg_{i_1}h^{-1},g_{i_1}^{-1}yg_{i_1}]\neq 1$. 
 
 We can thus suppose that $x$ and $y$ are not in the same conjugated component of $G$.
 We set $z=u^{nm}\in G(x)\cap G(y)$. Remark that the corollary 4.1.4 of [4] implies that $z$ is not the trivial element since $x$ and $y$ are not in the same conjugate of $G$. Since $z=x^m=y^n$ and $[G(x),G(y)]=1$, we deduce that $[G(z),G(z)]=1$.

 Suppose that the reduced form of $z$ starts and ends with  a generator of $F[X_i,i\in I]$, then if $g$ is an element of $G$ distinct of the identity,  $[z,gzg^{-1}]\neq 1$. This contradicts the fact that $[G(z),G(z)]=1$. If the last letter of the reduced sequence of $z$ is an element of $G$, we can conjugate $z$ by an element of $h$ of $G$ such that the last letter of the reduced element $hzh^{-1}$ is  a  generator of $F[X_i, i\in I]$. We can thus suppose that the reduced sequence of $z$ is of the form $g_1X^{l_1}_{i_1}..g_nX^{l_n}_{i_n}, g_1\neq 1$.

 Let $g\neq 1$ such that $g g_1\neq 1$, 
 
Since $g_1X^{l_1}_{i_1}..g_nX^{l_n}_{i_n}$ and $gg_1X^{l_1}_{i_1}..g_nX^{l_n}_{i_n}g^{-1}$ commute and are not in the same conjugated component of $G$, there exists $w$ such that
$g_1X^{l_1}_{i_1}..g_nX^{l_n}_{i_n}=w^a$ and $gg_1X^{l_1}_{i_1}..g_nX^{l_n}_{i_n}g^{-1}=w^b$.
This is impossible since  if the reduced form of $w$ has at its ends two elements of $G$, the reduced form of $w^c, c\in Z$ has at is ends two elements of $G$, and if one end of $w$ is a power of $X_i,i\in I$, one end of $w^c$ is also a power of $X_i$.

\end{proof}

{\bf Corollary 2.1}

{\it Suppose that $G$ is not trivial and does not have $T_1$-divisors of zero, then $G[X_i,i\in I]$ is an $T_1$-integral domain.}

\bigskip

We show now the following proposition:

\bigskip
{\bf Proposition 2.1.}

{\it Let $I$ and $J$ be two normal subgroups of $(H,f_G^H)$, we have $V^1_H([I,J])=V^1_H(I)\bigcup V^1_H(J)$ and $V^2_H(I\cap J)=V^2_H(I)\bigcup V^2_H(J)$.

Let $E(I_a)_{a\in A}$ be the normal subgroup of $H$ generated by the family of normal subgroups $(I_a)_{a\in A}$   of $H$, $V^i_H(E(I_a)_{a\in A})=\cap_{a\in A}V^i_H(I_a)$ $i=1,2$.}

\begin{proof}

Firstly, we show that $V^1_H([I,J])=V ^1_H(I)\bigcup V^1_H(J)$.
Let $P$ be an element of $V^1_H(I)\bigcup V^1_H(J)$, the fact that $P$ is an element of $V^1_H(I)$ or $V^1_H(J)$ implies that $P$ contains $I$ or $J$, henceforth $P$ contains $[I,J]$. On the other hand, let $P$ be an element of $V^1_H([I,J])$. Suppose that $P$ does not contain $I$ and $J$. This implies that there exist $x\in I$, and $y\in J$ which are not elements of $P$. Since $I$ and $J$ are normal subgroups of $H$, $G(x)\subset I$, and $G(y)\subset J$. This implies that $[G(x),G(y)]\subset [I,J]\subset P$. Since $P$ is a $T_1$-prime ideal, we deduce that $x\in P$ or $y\in P$. This is a contradiction. The proof that $V^2_H(I\cap J)=V^2_H(I)\bigcup V^2_H(J)$ is similar to the previous proof.

Now we show that $V^i_H(E(I_a)_{a\in A})=\cap_{a\in A}V^i_H(I_a)$. Let $P$ be an element of $V^i(E(I_a)_{a\in A})$, for every $a\in A$, $I_a\subset E(I_a)_{a\in A}\subset P$. This implies that $P\in \cap_{a\in A}V^i_H(I_a)$. Conversely, suppose that $P$ is an element of $\cap_{a\in A}V^i_H(I_a)$. This implies that for every $a\in A$, $I_a\subset P$ and henceforth, $E(I_a)_{a\in A}\subset P$.
\end{proof}

Remark that $V^i_H(1)=Spec^i_G(H)$ and $V^i_H(H)$ is empty. This implies that the family of subsets $V ^i_H(I)$, where $I$ is a normal subgroup of $H$ is the family of closed subsets of a topology defined on $Spec^i_G(H)$, $i=1,2$ that we call the $T_i$-topology.  Remark also that the proposition 2.1 implies that $V^1_H([H,H])=V^1_H(H)\bigcup V^1_H(H)$ is empty.

\medskip

We also have:

\bigskip

{\bf Proposition 2.2.}
{\it Let $(H,f_G^H)$ be an object of $C(G)$, for every normal subgroups $I,J$ of $H$, $V^1_H(I\cap J)=V^1_H([I,J])$
}

\begin{proof}
Since $[I,J]\subset I\cap J$, we deduce that $V^1_H([I\cap J)\subset V^1_H([I,J])$. On the other hand, Let $P$ be an element of $V^1_H([I,J])$. Let $x$ be an element of  $I\cap J$. We have $G(x)\subset I\cap J$ implies that $[G(x),G(x)]\subset [I,J]\subset P$, since $P$ is prime, we deduce that $x\in P$.
\end{proof}

{\bf Proposition 2.3.}

{\it Let $h\in H$, we denote by $V^i_h$ the set of $T_i$-prime ideals which contain $h$ and by $U^i_h$ its complementary subset $i=1,2$. The set $\{ U^i_h, h\in H\}$ is a base of the $T_i$-topology.
}

\begin{proof}
Let $n(h)$ be the normal subgroup generated by $h$, $V^i_h=V_H(n(h))$. This implies that $V^i_h$ is closed and $U^i_h$ is open. Let $V^i_H(I)$ be a closed subset of $Spec^i_G(H)$ and $U_I$ its complementary subset. Let $P\in U^i_I$, there exists an element $h\in I$ which is not an element of $P$. This implies that $P\in U^i_h$. Let $Q$ be any element of $U^i_h$. Since $Q$ does not contain $h$, we deduce that it does not contain $I$. This implies that $U^i_h\subset U^i_I$ and $U^i_I=\bigcup_{h\in I}U^i_h$.
\end{proof}

\bigskip

{\bf Proposition 2.4.}

{\it Let $H$ be a $G$-group and $I_1,..,I_n$ be normal subgroups of $H$. Suppose that  $I_1\cap I_2..\cap I_n=1$, and for $j\neq k$, $I_jI_k=H$,  then $Spec^i_G(H)$ is the disjoint union of the elements of the set $\{V^i_H(I_i),=1,...,n\}, i=1,2$.}

\begin{proof}
The fact that $I_i\cap I_2...\cap I_n=1$ implies that $V^i_H(I_1\cap I_2\cap..\cap I_n)=V^i_H(I_1)\bigcup..\bigcup V^i_H(I_n)=V^i_H(1)=Spec^i_G(H)$. We deduce from the fact that $I_jI_k=H$ if $j\neq k$, that $V^i_H(I_jI_k)=V^i_H(I_j)\cap V^i_H(I_k)=V^i_H(H)$. This implies that $V^i_H(I_jI_k)$ is empty since $V_H(H)$ is empty.
\end{proof}

\bigskip
{\bf Definition 2.2.}

{\it Let $U$ be an open subspace of $Spec^i_G(H)$. The $T_i$-radical $rad^i_H(U)$  is the intersection of the elements of $U, i=1,2$. We will often denote $rad^i_H(Spec^i_G(H))$ by $rad^i_G(H)$. 

Let $P$ be an element of $Spec^i_G(H)$, $rad_G^i(P)$ the $T_i$-radical of $P$ is the subgroup of $H$ generated by $\{rad^i_G(U)$, $U$ is an open subset such that $P\in U\}$. }

\medskip

Remark that
the adherence of $U$ for the $T_i$-topology is $V_H(rad^i_H(U))$.

\bigskip

{\bf Proposition 2.5.}

{\it Let $I$ be a normal subgroup, suppose that $\bigcap_{P\in V^i_H(I)}P=I$, then $V^i_H(I)$ is irreducible if and only if $I$ is a prime ideal $i=1,2$. In particular, $Spec^i_G(H)$ is irreducible if and only if $rad^i_G(H)$ is a prime ideal.}

\begin{proof}
We are going to proof this result for the $T_1$-topology, the corresponding proof for the $T_2$-topology is similar. Suppose that $V^1_H(I)$ is irreducible and $I$ is not a prime ideal; there exist elements $x,y\in H$ such that $[G(x),G(y)]\subset I$ and neither $x$ nor $y$ is an element of $I$. Let $P$ be an element of $V^1_H(I)$, $[G(x),G(y)]\subset I\subset P$. Since $P$ is a $T_1$-prime ideal, $x\in P$ or $y\in P$. Consider $V^1_x$ (resp. $V^1_y$) the subset of $T_1$-prime ideals containing the normal subgroup of $H$ generated by $x$ (resp. $y$). 
The closed subspace $V^1_H(I)$ is the union of $V^1_x\cap V^1_H(I)$ and $V^1_y\cap V^1_H(I)$, since it is irreducible, we deduce that $V^1_x\cap V^1_H(I)=V^1_H(I)$ or $V^1_y\cap V^1_H(I)=V^1_H(I)$. If $V^1_x\cap V^1_H(I)=V^1_H(I)$, $x\in P$ for every $P\in V^1_H(I)$. This implies that $x\in \cap_{P\in V^1_H(I)}P=I$. Contradiction.

Conversely, suppose that $I$ is a $T_1$-prime ideal and $V^1_H(I)$ is the union of the closed subsets $V^1_H(I)\cap V^1_H(J)$ and $V^1_H(I)\cap V^1_H(K)$. The $T_1$-prime ideal $I$ is an element of $V^1_H(I)\cap V^1_H(J)$ or $V^1_H(I)\cap V^1_H(K)$. Suppose $I$ is an element of $V^1_H(I)\cap V^1_H(J)$. This implies that $J\subset I$. We deduce that for every element $P\in V^1_H(I), J\subset I\subset P$. This implies $P\in V^1_H(J)$ and $V^1_H(I)\cap V^1_H(J)=V^1_H(I)$.
\end{proof}

Recall that a topological space is the union of its irreducible components. We can write $Spec^i_G(H)=\bigcup_{j\in J} V^i_H(P_j)$, where $V^i_H(P_j)$ is an irreducible component of $Spec^i_G(H), i=1,2$. Moreover, If we assume that $P_i=\cap_{Q\in V^i_H(P_i)}Q$. The proposition 2.5 implies that $P_i$ is a prime. We have the following result:

\bigskip
{\bf Proposition 2.6.}

{\it Let $H$ be an element of $C(G)$ and $V^i_H(P_j)$ an irreducible component of $Spec^i_G(H), i=1,2$ such that $P_j=\cap_{Q\in V_H(P_j)}Q$. For every element $Q$  of $V_H(P_j)$, every open subset  which contains $Q$ contains $P_j$.}

\begin{proof}
Suppose that $U$ is an open subset of $Spec^i_G(H)$ which contains $Q$ and $P_j$ is not an element of $U$; then $P_j$ is an element of the complementary space of $U$ which is a closed subset $V^i_H(I)$.This implies that $I\subset P_j\subset Q$. We deduce that $Q\in V^i_H(I)$. Contradiction.
\end{proof}

We say that the object $H$ of $C(G)$ is $T_i$-irreducible $i=1,2$ if and only if its radical is a $T_i$-prime ideal. The proposition 2.6 implies:

\bigskip
{\bf Corollary  2.2.}

{\it Suppose that $H$ is irreducible, $rad^i_G(H)$ is dense in $Spec^i_G(H)$, and  for every element $P$ of $Spec^i_G(H), rad^i_G(P)=rad^i_G(H), i=1,2$.}

\begin{proof}
Suppose that $H$ is irreducible, since $rad^i_G(H)=\bigcap_{P\in Spec^i_G(H)}P$, the proposition 2.5 implies that $rad^i_G(H)$ is a prime ideal, and the proposition 2.6 implies that every non empty open subset of $Spec^i_G(H)$ contains $rad^i_G(H)$. Henceforth the only non empty closed subset which contains $rad^i_G(H)$ is $Spec^i_G(H)$. This is equivalent to saying that $rad^i_G(H)$ is dense in $Spec^i_G(H)$.  Since $rad^i_G(H)$  is contained in every non empty open subset $U$ of $Spec^i_G(H)$, $rad^i_G(U)=rad^i_G(H)$, and henceforth, $rad^i_G(P)=rad^i_G(H)$.
\end{proof}

Let $H$ be an irreducible object of $C(G)$, the element $rad^i_G(H)$ is called the generic element of $Spec^i_G(H)$.

\bigskip
\bigskip

{\bf 3. $G$-varieties.}

\bigskip

Recall that a space $X$ is Noetherian if and only if for every sequence $(F_n)_{n\in \mathbb{N}}$ of closed subsets of $X$, such that $F_n\supset F_{n+1}$, there exists an integer $M$ such that $F_n=F_M$ for every $n\geq M$. 

\medskip

{\bf Definition 3.1.}

{\it A $G$-group $H$ is finitely $G$-generated if and only if  there exist elements $h_1,...,h_n\in H$ such that $G,h_1,..,h_n$ generate $H$. In this case, we say that $Spec^i_G(H)$ is an affine $G$-scheme $i=1,2$.}

\bigskip

The scheme $Spec^i_G(H)$ is an affine $G$-scheme if and only if there exists  a surjective $G$-morphism $i_H:G[X_1,..,X_n]\rightarrow H$. We denote by $I_H$ the kernel of $i_H$.

Let $G^n$ be the product of $n$-copies of $G$, each element $x=(g_1,..,g_n)$ of $G^n$ defines a morphism $e_x:G[X_1,..,X_n]\rightarrow G$ by $e_x(g)=g, g\in G,  e_x(X_i)=g_i$. Let $P$ be an element of $G[X_1,..,X_n]$, we denote by $f_P:G^n\rightarrow G$ the function defined by $f_P(x)=e_x(P)$. Such a function will be called a polynomial function. 
Let $g$ be an element of $G$, we denote by $c_g$ the polynomial function  such that for every $x\in G^n$, $c_g(x)=g$.

\bigskip

{\bf Definition 3.2.}

{\it Let  $I$ be a normal subgroup of $G[X_1,..,X_n]$, the affine $G$-variety  associated to the affine spectrum $G[X_1,..,X_n]/I$ is the subset $Var_G(I)=\{x\in G^n: $  for every $f\in I,  f(x)=1\}$.

Let $P$ be an element of $G[X_1,..,X_n]$, a zero of $P$ is an element $x\in G^n$ such that $f_P(x)=1$.  For every element $x$ of $G^n$, we denote by $I_x=\{P\in G[X_1,..,X_n]: f_P(x)=1\}$.
}

\bigskip
{\bf Proposition 3.1.}

{\it The ideal $I_x$ is maximal if and only if $G$ is simple.}

\begin{proof}
Suppose that $I_x$ is maximal, let $N$ be a normal subgroup of $G$, the set defined by $I_{N,x}=\{P\in G[X_1,..,X_n], f_P(x)\in N\}$ is a normal subgroup containing $I_x$, since $I_x$ is maximal, $I_{N,x}=I_x$ or $I_{N,x}=G[X_1,..,X_n]$. Suppose that $I_{N,x}=I_x$. This implies that for every $n\in N, f_{c_n}\in I_x$, this implies that $c_n(x)=n=1$ and $N=\{1\}$. If $I_{N,x}=G[X_1,..X_n]$, for every $g\in G, c_g\in I_{N,x}$. This implies that $c_g(x)=g\in N$. We deduce that $N=G$. 

Conversely, suppose that $G$ is simple, let $J$ be an ideal containing $I_x$, $J_x=\{f_P(x),P\in J\}$ is a normal subgroup of $G$. Since $G$ is simple, it implies that $J_x=\{1\}$ or $J_x=G$. If $J_x=\{1\}$, then $J_x=I_x$. If $J_x=G$. Let $P\in G[X_1,..,X_n]$ there exists $Q\in J$ such that $f_P(x)=f_Q(x)$. We have $P=PQ^{-1}Q$. Since $f_{PQ^{-1}}(x)=1$, we deduce that $PQ^{-1}\in I_x$ and henceforth $PQ^{-1}Q\in J$. This implies that $J=G[X_1,..,X_n]$.
\end{proof}

We also  have:

\bigskip
{\bf Proposition 3.2.}

{\it Suppose that $G$ is $T_i$-integral, then for every $x\in G^n$, $I_x$ is a $T_i$ prime ideal, $i=1,2$.}

\begin{proof}
We give the proof if $G$ is $T_1$-integral. The proof when $G$ is $T_2$-integral is obtained mutatis mutandis.
Let $f,f'\in G[X_1,...,X_n]$, suppose that $[G(f),G(f')]\subset I_x$. This implies that  $[G(f(x)),G((f'(x))]=1$. The hypothesis of the proposition implies that $f(x)=1$ or $f'(x)=1$. We deduce that $f\in I_x$ or $f'\in I_x$.
\end{proof}

The following proposition can be compared to [1] lemma 6. 

\bigskip
{\bf Proposition 3.3.}

{\it Suppose that $G$ is $T_1$-integral (resp. $T_2$-integral), then for every normal subgroups $I,J$ of $G[X_1,...,X_n]$, $Var_G(I)\bigcup Var_G(J)=Var_G([I,J])$ (resp. $Var_G(I)\bigcup Var_G(J)=Var_G(I\cap J)$).

For every family of normal subgroups $(I_a)_{a\in A}$ of $G[X_1,...,X_n]$, $\cap_{a\in A}Var_G(I_a)=Var(gr(I_a)_{a\in A})$ where $gr(I_a,a\in A)$ is the normal sugroup of $G$ generated by $(I_a)_{a\in A}$}

\begin{proof}
We are going to proof the result if $G$ is $T_1$-integral. The proof if $G_2$ is $T_2$-integral is obtained mutatis mutandis.

Let $x\in Var_G(I)\bigcup Var_G(J)$, and $f\in [I,J]$, if $x\in Var_G(I)$, $f\in [I,J]\subset I$, we deduce that $f(x)=1$. A similar argument shows that if $x\in Var_G(J), f(x)=1$. We deduce that $Var_G(I)\bigcup Var_G(J)\subset Var_G([I,J])$. Let $x\in Var_G([I,J])$, suppose that $x$ is not an element of $Var_G(I)\bigcup Var_G(J)$. This implies that there exist $f\in I, f'\in J$ such that $f(x)$ and $f'(x)$ are distinct of $1$; $[G(f),G(f')]\subset [I,J]$, this implies that $[G(f(x)),G(f'(x))]=1$, since $G$ is integral, we deduce that $G(f(x))=1$ or $G(f'(x))=1$. Contradiction.

Let $x$ be an element of $\cap_{a\in A}Var_G(I_a)$, for every $a\in A, f\in I_a, f(x)=1$. This implies that $x\in Var(gr((I_a)_{a\in A})$. Conversely, let $x$ be an element of $Var(gr((I_a)_{a\in A})$, for every $a\in A, f\in I_a$, $f(x)=1$ since $I_a\subset gr((I_a)_{a\in A})$, we deduce that $Var_G(gr(I_a)_{a\in A})\subset \cap_{a\in A}Var_G(I_a)$. 
\end{proof}

Remark that the empty set is $Var_G(G[X_1,..,X_n])$. If $G$ is a $T_1$-integral domain (resp. $T_2$-integral domain), we can define on $G^n$ the topology whose closed subsets are $Var_G(I)$. In this case,  we can endow  an affine $G$-variety $V=Var_G(I), I\subset G[X_1,...,X_n]$ with the topology induced by the topology of $G^n$.

\bigskip
{\bf Definition 3.3.}

{\it Let $V$ be an affine $G$-variety defined by the ideal $I_V\subset G[X_1,..,X_n]$ and $U$ an open subset of $V$.  A function $f:U\rightarrow G$ is regular if and only if there exists $P\in G[X_1,..,X_n]$ such that for every $x\in U, f(x)=f_P(x)$. We denote by $O_V(U)$ the set of regular functions defined on $U$.

A morphism  between the affine $G$-varieties $V$ and $V'$ is a continuous map $\phi:V\rightarrow V'$ such that for  every regular function $f\in O_{V'}(V')$, the function $f\circ \phi:V\rightarrow G$ is a regular function.
}

Remark that if $V=Spec_G(G[X_1,..,X_n]/I_V)$ is an affine $G$-variety, then $O_V(V)=G[X_1,..,X_n]/I_V$. We denote by $p_V:G[X_1,..,X_n]\rightarrow O_V(V)$ the quotient morphism.

\bigskip

{\bf Proposition 3.4.}

{\it Let $V$ and $V'$ be two affine $G$-varieties, there exists a natural bijection between   $Hom_G(V,V')$, the set of morphisms of affine $G$-varieties between $V$ and $V'$, and $Hom_G(O_{V'}(V'),O_V(V))$ the set of morphisms of $G$-groups between $O_{V'}(V')$ and $O_V(V)$.}

\begin{proof}
Let $\phi:V\rightarrow V'$ be a morphism of $G$-varieties, we denote by $a(\phi):O_{V'}(V')\rightarrow O_V(V)$ the morphism of $G$-groups defined by $a(\phi)(f)=f\circ \phi$.

Let $\psi:O_{V'}(V')\rightarrow O_V(V)$ be a morphism of $G$-groups. Suppose that the group of regular functions of $V$ (resp. $V'$) is $G[X_1,..,X_n]/I_V$  $(resp. G[X_1,..,X_{n'}]/I_{V'})$, we denote by $[X_i]_Y$ the image of $X_i$ by the quotient map $G[X_1,..,X_{n'}]\rightarrow G[X_1,..,X_{n'}]/I_{V'}$. For every $x\in V$, we write $b(\psi)(x)=(\psi([X_1]_Y)(x),..,\psi([X_{n'}]_Y)(x))$.  Let $P\in I_{V'}$, $f_P(b(\psi(x))=f_{\psi(P)}(x)=1$. We deduce  that the image of  $b(\psi)$ is contained in $V'$. Let $I$ be an ideal of $G[X_1,..,X_{n'}]$, $b(\psi)^{-1}(Var_{G}(I))=Var_G(\psi(I))$. This implies that $b(\psi)$ is continue.

The morphism $a:Hom_G(V,V')\rightarrow Hom_G(O_{V'}(V'),O_V(V))$ is the inverse of the morphism $b:Hom_G(O_{V'}(V'),O_V(V))\rightarrow Hom_G(V,V')$. To see this, remark that $a(b(\psi)[X_i]_Y(x))=[X_i]_Y(b(\psi)(x))=\psi([X_i]_Y)(x)$. On the other hand, let $b(a(f))(x)=(([X_1]_Y(f(x)),..,([X_{n'}]_Y(f(x))=f(x).$
\end{proof}

\medskip

{\bf Proposition 3.5.}

{\it Let $G$ be a $T_i$-integral domain $i=1,2$, and $V$ an irreducible affine $G$-variety. We denote by $I_V=\cap_{x\in V}I_x$.  The ideal $I_V$ is a $T_i$-prime ideal.}

\begin{proof}

We are going to prove this result for the $T_1$-topology, the proof for the $T_2$-topology is obtained mutatis mutandis.

Let $f,f'$ be elements of $G[X_1,...,X_n]$ such that $[G(f),G(f')]\subset I_V$. Consider an element $x$ of $V$. The ideal $I_x$ contains $I_V$ and henceforth $[G(f),G(f')]$. The proposition 3.2 implies that $f\in I_x$ or $f'\in I_x$. Let $n(f)$ (resp. $n(f')$) the normal subgroup of $G[X_1,...,X_n]$ generated by $f$ (resp. $f'$); $V=(V\cap V_G(n(f)))\bigcup (V\cap V_G(n(f')))$. Since $V$ is irreducible, we deduce that either $V\cap V_G(n(f))=V$ or $V\cap V_G(n(f'))=V$. If $V=V\cap V_G(n(f))$, we deduce that $V\subset V(n(f))$, for every $x\in V, f\in I_x$, henceforth, $f\in\cap_{x\in V}I_x=I_V$. If $V=V\cap V_G(n(f'))$, we deduce that $V\subset V(n(f))$, for every $x\in V, f'\in I_x$, henceforth, $f'\in\cap_{x\in V}I_x=I_V$. This is equivalent to saying that $I_V$ is a prime ideal.
\end{proof}

\bigskip
\bigskip

{\bf 4. The structural sheaf.}

\bigskip

Let $(H,f_H)$ be an object of $C(G)$, and  $U$ an open subset of $Spec^i_G(H), i=1,2$.
 The structural sheaf $O^i_H$ of $Spec^i_G(H)$ is defined as follows: for every open subset $U$ of $Spec^i_G(H)$, $O^i_H(U)$ is a subset of the set of functions $U\rightarrow \bigcup_{P\in U}H/P$ such that:
 
For every $s\in O^i_H(U), P\in U, s(P)\in H/P$,
 
  There exists an open cover $(U_j)_{j\in J}$ of $U$, an element $h_j\in H$ such that for every $P\in U_j, s(h)=l_P(h_j)$.

\bigskip
{\bf Proposition 4.1.}

{\it Let $P$ be an element  of $Spec^i_G(H), i=1,2$, the stalk ${O_H}^i_P$ of $O^i_H$ at $P$ is a quotient of  $H/rad^i_H(P)$.
}

\begin{proof}

Let $q_P:H\rightarrow {O_H}^i_P$ be the morphism of groups such that for every element  $h$  of $H$,  $q_P(h)$ is the image of $p_{Spec^i_G(H)}(h)$ in ${O_H}^i_P$ where $p_{Spec^i_G(H)}:\rightarrow H/rad^i_G(H)$ is the quotient map. We are going to show that $q_P$  is surjective and its kernel is $rad^i_H(P)$. 

For every $g\in {O_H}^i_P$, there exists an open subset $U$ containing $P$, $ h\in H$ such that $g$ is the image of $r_{U,Spec^i_G(H)}(h)$ in ${O_H}^i_P$. This implies that $g=q_P(h)$. We deduce that $q_P$ is surjective.

Let us show now that the kernel of $q_P$ is contained in $rad^i_H(P)$. Suppose that $q_P(h)=1$. This implies that there exists an open subset $P\in U$ such that $r_{U,Spec^i_G(H)}(h)=1$, this implies that $h\in rad^i_H(U)\subset  rad^i_H(P)$.
\end{proof}

\bigskip

{\bf Theorem 4.1.}

{\it Let $H$ be an element of $C(G)$, write $Spec^i_G(H)=\bigcup_{j\in J}V_H(P_j)$ where $V_H(P_j)$ is an irreducible component of $Spec^i_G(H), i=1,2$. Suppose that $\cap_{j\in J}V_H(P_j)$ is not empty, then   $O^i_H(Spec^i_G(H))=H/rad^i_G(H)$.
}

\begin{proof}

We are going to show that $p_{Spec^i_G(H)}$ is surjective; since its kernel  is $rad^i_G(H)$, this will show that $O^i_H(Spec^i_G(H))=H/rad^i_G(H)$.

Let   $P$ be an element of $\cap_{j\in J}V^i_H(P_j)$ and $h\in O^i_H(Spec^i_G(H))$, there exits an open subset $U$ of $Spec^i_G(H)$ containing $P$, an element $h_U\in H$ such that for every element $Q\in U$, $h(Q)=l_Q(h_U)$. Since $P\in V^i_H(P_j)$ for every $j\in J$, the proposition 2.6 implies that $P_j\in U$. We deduce that $h(P_j)=l_{P_j}(h_U)$. Let $Q'$ be any element of $Spec^i_G(H)$. There exits an open subset $U'$ of $Spec^i_G(H)$ which contains $Q'$, an element $h_{Q'}\in H$ such that for every element $Q\in U'$, $h(Q)=l_Q(h_{Q'})$.  There exists $j_{Q'}\in J$ such that $Q'\in V_H(P_{j_{Q'}})$. The proposition 2.6 implies that $P_{i_{Q'}}\in U'$. We deduce that $h(P_{i_{Q'}})=l_{P_{i_{Q'}}}(h_{Q'})=l_{P_{i_{Q'}}}(h_U)$. Since $h(Q')$ is the image of $h(P_{j_{Q'}})$ by the canonical projection $H/P_{j_{Q'}}\rightarrow H/Q'$, we deduce that $h(Q')=l_{Q'}(h_U)$. This implies that $p_{Spec^i_G(H)}$ is surjective.
\end{proof}

\bigskip
{\bf Corollary 4.1.}

{\it Suppose that $Spec^i_G(H), i=1,2$ is irreducible, then $O^i_H$ is a simple sheaf and $O_H(Spec^i_G(H))=H/rad^i_G(H)$.}

\begin{proof}
The theorem 4.1 implies that $O_H(Spec^i_G(H))=H/rad^i_G(H)$ if $Spec_G(H)$ is irreducible. It remains to show that for every open subset $U$ of $Spec^i_G(H)$, $O^i_H(U)=H/rad^i_G(H)$ and the restriction morphisms are the identities. Let $h\in O^i_H(U)$, for every $P\in U$,  there exists an open subset $V\subset U$ containing $P$, $h'\in H$ such that for every $Q\in V$, $h(Q)=l_{Q}(h')$. In particular, $h(rad^i_G(H))=l_{rad^i_G(H)}(h')$, we deduce that $h(P)$ is the image of $h(rad^i_G(H))$ by the quotient map $H/rad^i_G(H)\rightarrow H/P$. This implies that $h=r_{U,Spec^i_G(H)}(h(rad^i_G(H))$. We deduce that $r_{U,Spec^i_G(H)}$ is surjective; $r_{U,Spec^i_G(H)}(h)=1$ implies that $r_{U,rad^i_G(H)}(h)(rad^i_G(H))=h=1$. This implies that $r_{U,Spec^i_G(H)}$ is injective, and henceforth is an isomorphism.
\end{proof}

\bigskip
\bigskip

{\bf 5. $G$-schemes.}

\bigskip

In this section, we are going to study some properties of $T_i$-$G$-schemes, $i=1,2$.

\bigskip
{\bf Definition 5.1.}

{\it A $G$-space $(X,O_X)$, is a topological space $X$ endowed with a sheaf $O_X$, such that for every open subset $U$ of $X$, $O_X(U)$ is an object of $C(G)$, and   the restriction morphisms of the sheaf $O_X$ are morphisms of $G$-groups.

A local $G$-space is defined by a $G$-space $(X,O_X)$, and for every element $x\in X$, a prime ideal $P_x$ of ${O_X}_x$, the stalk of $x$.

Let $(X,O_X)$ and $(Y,O_Y)$ be two $G$-spaces. A morphism $(f,f^\sharp):(X,O_X)\rightarrow (Y,O_Y)$ of $G$-spaces is defined by a continuous map $f:X\rightarrow Y$, and a morphism of sheaves $f^\sharp:O_Y\rightarrow f^\sharp O_X$, such that for every open subsets $U\subset V$ of $Y$, $f^{\sharp}(U)$ is a morphism of $G$-groups and the following square is commutative:

$$
\matrix{O_Y(V)&{\buildrel{f^\sharp(V)}\over{\longrightarrow}}& O_X(f^{-1}(V))\cr
r_{U,V}\downarrow &&\downarrow r_{f^{-1}(U),f^{-1}(V)}\cr
O_Y(U)&{\buildrel{f^\sharp(U)}\over{\longrightarrow}}& O_X(f^{-1}(U))}
$$

A morphism between two local $G$-spaces $(f,f^{\sharp}):(X,O_X)\rightarrow (Y,O_Y)$ is a morphism of $G$- spaces such that: for every $x\in X$, the inverse image of $P_x$ by the  morphism ${O_Y}_{f(x)}\rightarrow {O_X}_x$ between the stalks of $f(x)$ and $x$ induced by $f^{\sharp}$ is $P_{f(x)}$.
}

\bigskip

{\bf Definition 5.2.}

{\it The topological space $Spec^i_G(H), i=1,2$ endowed with the sheaf $O^i_H$ is a $G$-space also called an affine $G$-spectrum. It will be often  denoted by $(Spec^i_G(H),O^i_H)$. 

An open subset $U$ of an affine $G$-spectrum $(Spec^i_G(H),O^i_H)$ endowed with the restriction of $O^i_H$ is called a quasi-affine $G$-spectrum.
}

\bigskip

Remark that the group space $(Spec^i_G(H),O^i_H)$ is naturally endowed with the structure of  a local $G$-space, such that, for every $P\in Spec^i_G(H)$, the prime element of the stalk of $P$ associated to $P$ is  the image of $P$ by the quotient morphism $H\rightarrow {O^i_{Spec_G(H)}}_P$.

\bigskip

{\bf Proposition 5.1.}

{\it Let $U$ be an open subset of the affine $G$-spectrum $Spec^i_G(H),i=1,2$ and $P$ an element of $U$, the subgroup $I^i_P(U)=\{h\in O^i_H(U): h(P)=1\}$ is a prime ideal of $O^i_H(U)$.
}

\begin{proof}
We are going to proof this result for the $T_1$-topology, the argument for the $T_2$-topology is obtained mutatis mutandis.

Let $h_1,h_2\in O^i_H(U)$, suppose that $[G(h_1),G(h_2)]\in I_P$. There exist open subsets $U_j, j=1,2$ such that $U_j$ contains $P$ and there exists $e(h_j)\in H$ such that for every $Q\in U_j$, $h_j(Q)=l_Q(e(h_j))$. Since $P\in U_1\cap U_2$, we deduce that for every $g_1,g_2\in G$, $[g_1h_1g_1^{-1},g_2h_2g_2^{-1}](P)=l_P([g_1e(h_1)g_1^{-1},g_2e(h_2)g_2^{-1}])=1$. This implies that $[g_1e(h_1)g_1^{-1},g_2e(h_2)g_2^{-1}]$ is an element of $P$. Since $P$ is a prime ideal, we deduce that $e(h_1)\in P$ or $e(h_2)\in P$. This implies that $h_1\in I^i_P$ or $h_2\in I^i_P$.
\end{proof}

Recall that if $U$ is an open subset of $Spec^i_G(H)$, an element $h\in H/rad^i_G(U)$ induces an element $s(h)\in O_H(U)$ such that for every $Q\in U, s(h)(Q)=l_Q(h)$. The element $s(h)\in O^i_H(U)$ will be called a complete section.

\bigskip
{\bf Proposition 5.2.}

{\it Let $(Spec^i_G(H),O^i_H)$ and $(Spec^i_G(H'),O^i_{H'})$ be two affine $T_i$-$G$-spectra $i=1,2$, a $G$-morphism $f:H\rightarrow H'$ induces a morphism of local $G$-spaces $(i_f,i_f^{\sharp}):(Spec^i_G(H'),O^i_{H'})\rightarrow (Spec^i_G(H),O^i_H)$ defined as follows:

For every $P\in Spec^i_G(H')$,  $i_f(P)=f^{-1}(P)$.

Let $U$ be an open subset of $Spec^i_G(H)$ and $h$ an element of $O^i_H(U)$, for every $P\in U$, there exists an open subset $U_P$ containing $P$, an element $h_P\in H$ such that for every $Q\in U_P$, $h(Q)=l_Q(h_P)$.  For every $L\in {i_f}^{-1}(P), i_f^{\sharp}(h)(L)=l_L(f(h_P)).$ 
}

\begin{proof}
Firstly, we show that $i_f$ is well-defined. Let $P\in Spec^i_G(H')$,suppose that $f^{-1}(P)=H$, we deduce that the image of $f_{H'}=f\circ f_H$ is contained in $P$ and $G$ acts trivially on $H/P$. This is a contradiction with the fact that $H/P$ is not the trivial group and does not have divisors of zero. The morphism $f$ induces an injective  $G$-morphism $H/f^{-1}(P)\rightarrow H'/P$; since $H'/P$ is an integral domain, $H/f^{-1}(P)$ is also an integral domain; this implies that $f^{-1}(P)$ is an element of $Spec^i_G(H)$.
Let $C$ be a closed subset of $Spec^i_G(H)$, write $C=V^i_H(I)$ and $n(f(I))$ the normal subgroup of $H'$ generated by $f(I)$, ${i_f}^{-1}(V^i_H(I))=V^i_{H'}(n(f(I))$. This implies that $i_f$ is continuous.

Let $U$ be an open subset of $Spec^i_G(H)$,  $P$ an element of ${i_f}^{-1}(U)$ and $h$ an element of $O^i_H(U)$. Suppose that $f(P)$ is an element of the open subset $U_P$  (resp. $V_P$) such that there exists $h_P\in H$ (resp. $h'_P\in H$) such that for every $Q\in U_P$ (resp. $Q'\in V_p$) $h(Q)=l_Q(h_P)$ (resp. $h(Q')=l_{Q'}(h'_P)$). Since $O^i_H$ is a sheaf, for every $Q\in U_P\cap V_P$, $h(Q)=l_Q(h_P)=l_Q(h'_P)$. Let $Q\in f^{-1}(P)$, $l_Q(f(h_P))$ is the image of $l_P(h_P)$ by the morphism $H/P\rightarrow H'/Q$ induced by $f$. Since $l_P(h_P)=l_P(h'_P)$ we deduce that $i_f^{\sharp}(h)(Q)$ does not depend of the neighborhood $U_P$ of $P$ chosen and henceforth, is well defined.
\end{proof}

Conversely, we have the following:

\medskip

{\bf Proposition 5.3.}

{\it Let $H$ and $H'$ be irreducible objects of $C(G)$ such that $rad^i_G(H)=1$, and $rad^i_G(H')=1$, and $U$ an open subset of $Spec^i_G(H)$. A morphism of local  $G$-spaces $(f,f^{\sharp}):(U,{O^i_H}_{\mid U})\rightarrow (Spec^i_G(H'),O^i_{H'})$ is induced by a morphism of groups $H'\rightarrow H$.
}

\begin{proof}
We set $g=f^{\sharp}(Spec^i_G(H'))$. Since $Spec^i_G(H)$ and $Spec^i_G(H')$ are irreducible objects, the domain of $g$ is $H'$ and its codomain is $H$ since $rad_G^i(H)=1$ and $rad_G^i(H')=1$. We are going to show that $(f,f^{\sharp})$ is induced  by $g$. 

Let $P$ be an element of $U$, since $H$ and $H'$ are irreducible,  $O^i_H$ and $O^i_{H'}$ are simple sheaves (see corollary 4.1),  the stalk of $P$ is $H$, the stalk of $f(P)$ is $H'$, and the morphism induced by $f$ on these stalks is $g$.  We deduce that $g^{-1}(P)=f(P)$ since $(f,f^{\sharp})$ is a local morphism of $G$-group spaces.
\end{proof}

\medskip

{\bf Corollary 5.1.}

{\it Let $H$ be an object of $C(G)$, and $I$ a normal subgroup of $H$. The quotient map $p_I:H\rightarrow  H/I$ induces an embedding of $G$-schemes $(i_{p_I},{i_{i_P}}^{\sharp}): (Spec^i_G(H/I),O^i_{H/I})\rightarrow (Spec^i_G(H),O^i_H), i=1,2$. 
In particular, if $I=rad^i_G(H)$, we obtain an isomorphism $(i_{P_{rad^i_G(H)}},{i_{P_{rad^i_G(H)}}}^{\sharp}):
(Spec^i_G(H/rad^i_G(H)),O^i_{H/rad^i_G(H)})\rightarrow (Spec^i_G(H),O^i_H)$.
}

\begin{proof}
The proposition 5.2 implies that $p_I$ induces a morphism $(Spec^i_G(H/I),O^i_{H/I})\rightarrow (Spec^i_G(H),O^i_H)$ which is an embedding since $p_I$ is surjective. To simplify the notations we are going to replace $p_{P_{rad^i_G(H)}}$ by $p$.  If $Q$ is an element of $Spec^i_G(H)$, $p(Q)$ is a prime ideal of $H/rad^i_G(H)$ and $i_P(p(Q))={p}^{-1}(p(Q))=Q$.
This implies that $i_P$ is bijective. Let $N$ be a normal subgroup of $H$, $i_{P}^{-1}(V^i_H(N))=V^i_{H/rad^i_G(H)}(p(N))$.
 This implies that $i_P$ is continue. If $U$ is an open subset of $Spec^i_G(H)$, $O^i_H(U)=O^i_{H/rad^i_G(H)}(i_P(U))$ and $i_P^{\sharp}(U)$ is the identity morphism. We deduce that $(i_P,i_P^{\sharp})$ is an isomorphism.
\end{proof}

\medskip

Remark that $Spec^i_G(H/I)=V^i_H(I)$, the previous corollary  enables to endow $V^i_G(I)$ with the structure of a $T_i$-$G$-spectrum.

\medskip

{\bf Definition 5.3.} 

{\it A morphism of $G$-spaces $(Spec^i_G(H),O^i_H)\rightarrow Spec^i_G(H'),O^i_{H'}), i=1,2$ induced by a $G$-morphism $H'\rightarrow H$ is called an affine morphism of $G$-spectra.

}

\medskip

{\bf Definition 5.4}

{\it Let $G$ be a group, a $T_i$-$G$-scheme $i=1,2$, is a local $G$-space $(X,O_X)$, such that  there exists an open cover $(U_j)_{j\in J}$ of $X$ such that for each $j\in J$, there exists an affine $G$-spectrum $ (Spec_G(H_j),O_{H_j})$, and an isomorphism of group spaces $(g_j,{g_j}^{\sharp}):(U_j,{O_X}_{\mid U_j})\rightarrow (Spec_G(H_j),O_{H_j})$.

A morphism of $T_i$-$G$-schemes $(f,f^\sharp):(X,O_X)\rightarrow (Y,O_Y)$ is a morphism of local $G$-spaces. 
}

Let $(f,f^\sharp):(X,O_X)\rightarrow (Y,O_Y)$ be a morphism of $T_i$-$G$-schemes, we will call $f$ the geometric component and $f^{\sharp}$ the algebraic component.

Let $x$ be an element of the $G$-scheme $(X,O_X)$, a chart containing $x$ is an open subset $U$ of $X$, such that $x$ is an element of $U$, and the restriction of $O_X$ to $U$ endows $(U,{O_X}_{\mid U})$ with the structure of a local $G$-space isomorphic to  an affine $T_i$-$G$-spectrum.

\medskip
{\bf Definition 5.5.}

{\it Let $(X_j,O_{X_j}), j=1,2$ be two $T_i$-$G$-schemes.  Suppose that there exist open subsets $U_j\subset X_j$ and an isomorphism $(f,f^{\sharp}): (U_1, {O_{X_1}}_{\mid U_1})\rightarrow (U_2,{O_{X_2}}_{\mid U_2})$. The gluing of $(X_1,O_{X_1})$ and $(X_2,O_{X_2})$ along $U_1$ and $U_2$ is the $G$-scheme $(X,O_X)$ such that $X$ is the quotient space 
$X_1\bigcup X_2/x\simeq f(x)$ endowed with the quotient topology. We denote by   $j_{X_j}:X_j\rightarrow X, j=1,2$ the canonical embedding map. A subset $U$ of $X$ is open if and only if $j_{X_j}^{-1}(U), j=1,2$ is open.

An element of $O_X(U)$ is a couple $(s_1,s_2)$ such that $s_1\in O_{X_1}(j^{-1}_{X_1}(U)), s_2\in O_{X_2}(j^{-1}_{X_2}(U))$, and such that $f^{\sharp}(j_{X_2}^{-1}(U)\cap U_2)({s_2}_{\mid j_{X_2}^{-1}(U)\cap U_2})={s_1}_{ j_{X_1}^{-1}(U)\cap U_1}$.}

\bigskip
{\bf Definition 5.6}

{\it A $T_i$-$G$-scheme $(X,O_X)$ is locally irreducible if and only if there exists an open covering $(U_j)_{j\in J}$ of $X$, such that for every $j\in J$, there exists an irreducible $G$-group $H_j$,  and an isomorphism of $G$-spaces $(U_j,{O_X}_{\mid U_j})\rightarrow (H_j,{O_{H_j}})$.
}

\medskip

We have the following result.

\bigskip
{\bf Theorem 5.1}

{\it Let $(X,O_X)$ be a locally irreducible $T_i$-$G$-scheme $i=1,2$, and $H$ an irreducible object of $C(G)$. There exists a bijection between the set of morphisms of $T_i$-$G$-schemes $(f,f^{\sharp}):(X,O_X)\rightarrow (Spec^i_G(H),O^i_H)$ and  $Hom_{C(G)}(O^i_H(Spec^i_G(H)),O^i_X(X))$ defined by $\Phi(f,f^{\sharp})= f^{\sharp}(Spec^i_G(H))$.
}

\begin{proof}
We are going to construct the inverse $\Psi$ of $\Phi$. Let $v^{\sharp}:O^i_H(Spec^i_G(H))\rightarrow O_X(X)$ be a $G$-morphism.  For every element $x$ of $X$, let $(U_j,{O_X}_{\mid U_j})$ be a chart containing $x$ which is isomorphic to  the  affine $G$-spectrum $(Spec^i_G(H_j),O^i_{H_j})$. For every $Q\in U_j$, we denote by $I^i_Q(U_j)=\{h\in O_X(U_j), h(x)=1\}$.  Consider $v_{U_j}$ the composition of $v^{\sharp}$ and the restriction $O_X(X)\rightarrow O_X(U_j)$. The proposition 5.1 shows that $I^i_Q(U_j)$ is a prime ideal of ${O_X}(U_j)$. We deduce that $v_{U_j}^{-1}(I^i_Q(U_j))$ is a prime ideal of $O^i_H(Spec^i_G(H))$. The corollary 4.1 shows that $O^i_H(Spec^i_G(H))=H/rad_G^i(H)$. We denote by $f_j(Q)$ the subset of complete sections of $v_{U_j}^{-1}(I^i_Q(U_j))$. It  is prime ideal of $H/rad^i_G(H)$. 

Let $(U_k,{O_X}_{\mid U_k})$ be another chart containing $x$,  we have: $v^{-1}_{U_j}(I^i_Q(U_j))=v^{-1}_{U_j\cap U_k}(I^i_Q(U_j\cap U_k))=v^{-1}_{U_k}(I^i_Q(U_k)$.This implies that $f_j(Q)=f_k(Q)$. We can define $f(Q)=f_j(Q)$.

Let $I$ be a normal subgroup of $H$, for every element $Q\in U_j$, $f(Q)\in V^i_H(I)$ if and only if the set of complete sections of $v_{U_j}^{-1}(I^i_Q(U_j))$ contains $I$, this is equivalent to say that $I^i_Q(U_j)$ contains $v_{U_j}(I)$. This last fact is equivalent to say that $n((v_{U_j}(I)))$ the normal subgroup generated by $v_{U_j}(I))$ is contained in $Q$. We deduce that $f^{-1}(V^i_H(I))\cap U_j=V^i_{H_j}(n(v_{U_j}(I))$. This implies that $f$ is continue.

Let $W$ be an open subset of $Spec^i_G(H)$, the  morphism $v^{\sharp}$ induces a morphism $v_W:O^i_H(W)\rightarrow O^i_X(f^{-1}(W))$ defined as follows: let $h\in O^i_H(W)$ and $Q\in f^{-1}(W)$, there exists a neighbourhood $U_{f(Q)}$ of $f(Q)$, $h'\in H$ such that for every $P\in U_{f(Q)}, h(P)=l_P(h')$.  we set $v_W(h)(Q)=v^{\sharp}(h')(Q)$.  The morphism $f$ is the geometric component of a morphism of $G$-schemes whose algebraic part is defined by the family of morphisms $(v_W)$. 

Remark that $v_{Spec^i_G(H)}=v^{\sharp}$; we deduce that $\Phi(\Psi(v^{\sharp}))=v^{\sharp}$.

Let $(f,f^{\sharp}):(X,O_X)\rightarrow (Spec^i_G(H),O^i_H)$ be a morphism between $T_i$-$G$-schemes.
Consider an element $P$ of an open subset $U_j\subset X$ such that $(U_j,{O_X}_{\mid U_j})$ is isomorphic to an  irreducible $T_i$-$G$-spectrum $(H_j,O_{H_j})$. We denote by $v_{U_j}$ the composition of $f^{\sharp}(Spec^i_G(H))$ and the restriction map $O_X(X)\rightarrow O_X(U_j)$. The set of global sections of $v_{U_j}^{-1}(I^i_P(U_j))$ is a prime ideal $Q$ of $H$ characterized by the fact that an element $h\in H$ is an element of $Q$ if and only if $v_{U_j}(h)(P)=1$.  Since $H$ and $H_j$ are irreducible, we can deduce like in  the proposition 5.3  that $Q=f(P)$. This implies that  $\Psi(\Phi(f,f^{\sharp}))=(f,f^{\sharp})$.
\end{proof}

Let $Spec^i_G(H), i=1,2$ be a Noetherian affine $T_i$-$G$-scheme. We can write $Spec^i_G(H)=\bigcup_{j=1}^{j=n}V_H(P_j)$ where $P_j$ is a prime ideal. Denote by $p_j:H\rightarrow H/P_j$  and $p_{jk}:H/P_k\rightarrow H/P_iP_j$ the canonical projections.

We have the following result:

\bigskip
{\bf Theorem 5.2.}

{\it Let $Spec^i_G(H)$ be a Noetherian affine $T_i$-$G$-scheme, write $Spec^i_G(H)=\bigcup_{j=1}^{j=n}V^i_H(P_j)$, $O^i_H(Spec^i_G(H))$ is the set $(g_1,...,g_n)\in \Pi_{j=1}^{=jn}H/P_j$ such that $p_{jk}(g_k)=p_{kj}(g_j)$ if the intersection $V^i_H(P_j)\cap V^i_H(P_k)$ is not empty.
}

\begin{proof}
We denote by $L_H$ the set $(g_1,...,g_n)\in \Pi_{j=1}^{j=n}H/P_j$ such that $p_{jk}(g_k)=p_{kj}(g_j)$ if the intersection $V^i_H(P_j)\cap V^i_H(P_k)$ is not empty. Let $h\in O^i_H(Spec^i_G(H))$, for every $j=1,..,n$ there exists a neighbourhood $U_j$ containing $P_j$, an element $h'_j\in H$ such that for every $Q\in U_j$, $h(Q) =l_{Q}(h'_j)$. We write $h_j=h(P_j)$. 

Suppose that $V^i_H(P_j)\cap V^i_H(P_k)$ is not empty and denote by $Q$ one of its element. Let $U_Q$ be an open subset containing $Q$ such that there exists $h_Q\in H$ such that for every $Q'\in U_Q$, $h(Q)=l_{Q'}(h_Q)$. The complementary space of $U_Q$ is a closed subset $V^i_H(I_Q)$, suppose that $P_j\in V^i_H(I_Q)$ (resp. $P_k\in V^i_H(I_Q)).$ This implies that $I_Q\subset P_j\subset Q$ (resp. $I_Q\subset P_k\subset Q$) contradiction. We deduce that $P_j\in U_Q$, $h(P_j)=l_{P_j}(h_Q)=h_j$ and $h(P_k)=l_{P_k}(h_Q)=h_k$. We deduce that $g_{jk}(h_k)=g_{jk}(l_{P_k}(h_Q)))=g_{kj}(l_{P_j}(h_Q)))=g_{kj}(h_i)$. We can define  the morphism $\Phi:O_H(Spec^i_G(H))\rightarrow L_H$ which associates $(h_1,...,h_n)$ to $h$. We are going to show that $\Phi$ is an isomorphism.

Suppose that $\Phi(h)=1$. Consider an element $Q\in Spec^i_G(H)$, and an open subset $U_Q$ such that there exists $h_Q\in H$ such that for every $Q'\in U_Q$, $h(Q')=l_{Q'}(h_Q)$. There exists $j\in \{1,...,n\}$ such that $Q\in V^i_H(P_j)$ and $P_j\in U_Q$. We deduce that $h(P_j)=l_{P_j}(h_Q)=1$, since $h(Q)$ is the image of $h_j$ by the canonical projection $U/P_i\rightarrow U/Q$, we deduce that $h(Q)=1$.This implies that $\Phi$ is injective.

Let $(h_1,..,h_n)$ be elements of $L_H$, define the element $h\in O^i_H(Spec^i_G(H))$ has follows: $h(P_j)=h_j$. Let $Q\in Spec^i_G(H)$, $U_Q$ an open subset containing $Q$. There exists $j$ such that $Q\in V_H(P_j)$. We have seen that $P_j$ is an element of $U_Q$; we denote by $h(Q)$ the image of $h_j$ by the quotient map $H/P_j\rightarrow H/Q$.
The element $h(Q)$ of $H/Q$ does not depend of the choice of $P_j$ such that $Q\in V_H(P_j)$. If $Q\in V_H(P_k), k\in\{1,...,n\}$, the image of $h_k$ by the canonical projection $H/P_k\rightarrow H/Q$ is the image of $p_{jk}(h_j)=p_{kj}(h_k)$ by the canonical projection $H/P_jP_k\rightarrow H/Q$. We deduce that $h$ is well defined and henceforth, $\Phi$ is an isomorphism.
\end{proof}

\bigskip

{\bf References.}

\medskip

1.  Baumslag, G,  Miasnikov, A.  Remeslennikov, V.N. Algebraic geometry over groups I. Algebraic sets and ideal theory. J. Algebra. 1999, 219, 1679.

 \smallskip
 
2. Griffiths, P.,  Harris, J. (2014). Principles of algebraic geometry. John Wiley and Sons.
 \smallskip
 
3. A.Grothendieck, \'El'ements de g\'eom\'etrie alg\'ebrique I.Publications math\'ematiques de l'I.H.E.S 4, 5-228
 \smallskip
 
4. Magnus, Karass, Solitar.  Combinatorial group theory, Dover publication 1976.
 \smallskip
 
5. Tsemo, Aristide. (2013). Some Properties of Schemes in Groups Theory and Top Couples. International Journal of Algebra, 7(1), 25-48.
\smallskip

\end{document}